\documentstyle[epsf,12pt]{article}
\def\Bbb R{{\rm \bf R}}
\def\proclaim#1{\vskip2mm{\bf #1}\em}
\def\endproclaim{\em \vskip2mm}
\def\tag#1{\eqno(#1)}
\def\gathered{\begin{array}{c}}
\def\endgathered{\end{array}}
\def\text{\mbox}

\begin{document}

\title {The enclosure method using a single point on
the graph of the response operator for the Stokes system}
\author{Masaru IKEHATA\footnote{
Laboratory of Mathematics,
Graduate School of Advanced Science and Engineering,
Hiroshima University,
Higashihiroshima 739-8527, JAPAN. 
e-mail address: ikehataprobe@gmail.com}
\footnote{Emeritus Professor at Gunma University, Maebashi 371-8510, JAPAN}
}
\maketitle

\begin{abstract}
An inverse obstacle problem governed by the Stokes system in the time domain is considered.
Two types of extraction formulae about the geometry of an unknown obstacle are given by using 
the most recent version of the time domain enclosure method in a unified style. 
Each of the formulae employs only a single set of velocity and the Cauchy stress fields over a finite time interval
on the surface of the region where a viscous incompressible fluid occupies.

\noindent
AMS: 35R30, 76M21, 76D07

\end{abstract}


\section{Introduction}

{\it The time domain enclosure method} with a single input is one of analytical methods in inverse obstacle problems
firstly initiated in \cite{I4} and developed in the series of articles
\cite{IW00, IEO2, IEO3} for an obstacle embedded in the {\it whole space}, and \cite{Iwave, IE05, IE06} 
for an obstacle embedded in a {\it bounded domain}.
However the pursuit of the most recent version of the time domain enclosure method developed in \cite{Iwave, IE06} 
is just beginning.
There are a lot of important, typical and interesting
inverse obstacle problems governed by various types of partial differential equations in the time domain.
In this paper, we consider how this new method can be realized in an inverse obstacle problem governed by the Stokes system
in the time domain.
Furthermore, we add new knowledge which was not covered by the previous research \cite{MST} under the framework 
of the previous version of the enclosure method in the time domain \cite{IFR}.

First we formulate the problem.
Let $\Omega$ be a bounded domain of $\Bbb R^3$ with $C^2$-boundary
and $D$ a nonempty open subset of $\Omega$ with $C^2$-boundary such that
$\overline D\subset\Omega$ and $\Omega\setminus\overline D$ is connected.
We denote by $\mbox{\boldmath $n$}$ the unit outward normal to $\partial\Omega$.

Let $0<T<\infty$.
Given $\mbox{\boldmath $f$}$ on $\partial\Omega\times [0,\,T]$
let $\mbox{\boldmath $u$}=\mbox{\boldmath $u$}(x,t)$ and $p=p(x,t)$ solve
$$\left\{
\begin{array}{ll}
\displaystyle
\rho\frac{\partial\mbox{\boldmath $u$}}{\partial t}
-\text{div}\,\sigma(\mbox{\boldmath $u$},p)=\mbox{\boldmath $0$},
&
\displaystyle
(x,t)\in(\Omega\setminus\overline D)\times\,]0,\,T[,
\\
\\
\displaystyle
\nabla\cdot\mbox{\boldmath $u$}=0,
&
\displaystyle
(x,t)\in(\Omega\setminus\overline D)\times\,]0,\,T[,
\\
\\
\displaystyle
\mbox{\boldmath $u$}=\mbox{\boldmath $0$},
&
\displaystyle
(x,t)\in\partial D\times\,]0,\,T[,\\
\\
\displaystyle
\mbox{\boldmath $u$}=\mbox{\boldmath $f$},
&
\displaystyle
(x,t)\in\partial\Omega\times\,]0,\,T[,\\
\\
\displaystyle
\mbox{\boldmath $u$}(x,0)=\mbox{\boldmath $0$},
&
\displaystyle
x\in\Omega\setminus\overline D.
\end{array}
\right.
\tag {1.1}
$$
Here $\rho$ denotes the density of the fluid occupying $\Omega\setminus\overline D$.
The pair $(\mbox{\boldmath $u$}, p)$ is that of the velocity and pressure of the fluid and
$\sigma(\mbox{\boldmath $u$},p)$ denotes the Cauchy stress tensor 
$$\displaystyle
\sigma(\mbox{\boldmath $u$},p)=-pI_3+2\mu\,\text{Sym}\nabla\mbox{\boldmath $u$},
$$
where $\mu$ denotes the viscosity of the fluid occupying $\Omega\setminus\overline D$.
It is assumed that $\rho$ and $\mu$ are given by positive constants.

Note that to ensure the existence of the solution of (1.1) we always choose $\mbox{\boldmath $f$}$ satisfying
$$\begin{array}{ll}
\displaystyle
\int_{\partial\Omega}\mbox{\boldmath $f$}(x,t)\cdot\mbox{\boldmath $n$}(x)\,dS=0,
& 0<t<T.
\end{array}
\tag {1.2}
$$

We consider

$\quad$

{\bf\noindent Problem.}  Fix $T$ and $\mbox{\boldmath $f$}$ (to be specified later).
Extract information about the geometry of $D$ from the pair 
of $\mbox{\boldmath $f$}$ and the corresponding the Cauchy force $\sigma(\mbox{\boldmath $u$},p)\mbox{\boldmath $n$}$
on $\partial\Omega$ over the time interval $]0,\,T[$.
In other words, {\it find} a suitable $\mbox{\boldmath $f$}$ in such a way that
the corresponding $\sigma(\mbox{\boldmath $u$},p)\mbox{\boldmath $n$}$ on $\partial\Omega$
over time interval $]0,\,T[$ yields some information
about the location and shape of $D$.

$\quad$

In this paper, we present two types of the solutions to this problem in a unified style.
It is an application of the idea introduced in \cite{Iwave} and the recent one in \cite{IE06}.

\subsection{Indicator function}

Let $\mbox{\boldmath $v$}=\mbox{\boldmath $v$}(x,t)$ and $q=q(x,t)$ satisfy
$$\left\{
\begin{array}{ll}
\displaystyle
\rho\frac{\partial\mbox{\boldmath $v$}}{\partial t}
-\text{div}\,\sigma(\mbox{\boldmath $v$},q)=\mbox{\boldmath $0$},
&
\displaystyle
(x,t)\in\Bbb R^3\times\,]0,\,T[,
\\
\\
\displaystyle
\nabla\cdot\mbox{\boldmath $v$}=0,
&
\displaystyle
(x,t)\in\Bbb R^3\times\,]0,\,T[.
\end{array}
\right.
\tag {1.3}
$$
Note that $\mbox{\boldmath $v$}$ satisfies
$$\displaystyle
\int_{\partial\Omega}\mbox{\boldmath $v$}\cdot\mbox{\boldmath $n$}\,dS
=\int_{\Omega}
\nabla\cdot
\mbox{\boldmath $v$}\,dx=0.
$$
Thus the vector field $\mbox{\boldmath $f$}$ given by
$$\begin{array}{ll}
\displaystyle
\mbox{\boldmath $f$}=\mbox{\boldmath $v$},
&
(x,t)\in\partial\Omega\times\,]0,\,T[,
\end{array}
\tag {1.4}
$$
satisfies (1.2).

$\quad$

{\bf\noindent Definition 1.1.} Let $(\mbox{\boldmath $u$},p)$ solve (1.1) with $\mbox{\boldmath $f$}$ given by (1.4).
Let $\tau>0$.
Define
$$\begin{array}{ll}
\displaystyle
I_T(\tau;\mbox{\boldmath $v$},q)
&
\displaystyle
=\int_{\partial\Omega}
\left(\sigma(\mbox{\boldmath $w$},\tilde{p})\mbox{\boldmath $n$}
-\sigma(\mbox{\boldmath $w$}_0,\tilde{q})\mbox{\boldmath $n$}\right)
\cdot\mbox{\boldmath $w$}_0\,dS,
\end{array}
\tag {1.5}
$$
where
$$\left\{
\begin{array}{ll}
\displaystyle
\mbox{\boldmath $w$}(x,\tau)=\int_0^Te^{-\tau T}\mbox{\boldmath $u$}(x,t)dt,
&
\displaystyle
\mbox{\boldmath $w$}_0(x,\tau)=\int_0^Te^{-\tau T}\mbox{\boldmath $v$}(x,t)dt,\\
\\
\displaystyle
\tilde{p}(x,\tau)=\int_0^Te^{-\tau T}p(x,t)dt,
&
\displaystyle
\tilde{q}(x,\tau)=\int_0^Te^{-\tau T}q(x,t)dt.
\end{array}
\right.
$$
We call the function $I_T(\tau;\mbox{\boldmath $v$},q)$ of independent variable $\tau$ defined by (1.5) the {\it indicator function}.

$\quad$

In principle, the indicator function can be computed from the pair $(\mbox{$f$}, \sigma(\mbox{\boldmath $u$},p)\mbox{\boldmath $n$})$
on $\partial\Omega$ over the time interval $]0,\,T[$ which is a point on the graph of the response operator
$$\begin{array}
{ll}
\displaystyle
\mbox{\boldmath $f$}\longmapsto\sigma(\mbox{\boldmath $u$},p)\mbox{\boldmath $n$}, & (x,t)\in\partial\Omega\times\,]0,\,T[.
\end{array}
$$

We show that, for suitable two choices of the pair $(\mbox{\boldmath $v$},q)$ with $q=0$ the asymptotic behaviour
of indicator function $I_T(\tau;\mbox{\boldmath $v$},q)$ as $\tau\rightarrow\infty$ yields some information about the geometry of obstacle $D$.
This is a solution to the problem mentioned above.

The point is: to generate the velocity and pressure of the fluid inside a given domain
we prescribe a special velocity field on the boundary of the domain given by solving
an evolution equation in the whole space.

The general outline of the method is as follows.

\noindent
(i)  Let $\Phi=\Phi(x,t)$ with $(x,t)\in\Bbb R^3\times\,[0,\,T]$ be a solution of the heat equation
$$\begin{array}{ll}
\displaystyle
\rho\frac{\partial\Phi}{\partial t}-\mu\Delta\Phi=0, & (x,t)\in\Bbb R^3\times\,]0,\,T[
\end{array}
\tag {1.6}
$$
such that $K=\text{supp}\,\Phi(\,\cdot\,,0)$ is compact and satisfies $K\cap\Omega=\emptyset$.

\noindent
(ii)  Define the vector field $\mbox{\boldmath $v$}$ by the formula
$$\begin{array}{ll}
\displaystyle
\mbox{\boldmath $v$}(x,t)=\nabla\times\{\Phi(x,t)\mbox{\boldmath $a$}\}, & (x,t)\in \Bbb R^3\times\,[0,\,T],
\end{array}
\tag {1.7}
$$
where $\mbox{\boldmath $a$}$ is an arbitrary constant unit vector.

We see that the pair of $\mbox{\boldmath $v$}$ and $q=0$ solves (1.3).

\noindent
(iii) Prescribe $\mbox{\boldmath $f$}$ given by (1.4) as the velocity field on $\partial\Omega$ 
over time interval $]0,\,T[$ and solve (1.1).

\noindent
(iv) Compute the indicator function $I_T(\tau;\mbox{\boldmath $v$},0)$.  Then, the asymptotic behaviour
of this indicator function as $\tau\rightarrow\infty$ yields the quantity $\text{dis}\,(D,K)$ for special $K$
which depends on the initial data $\Phi(\,\cdot\,,0)$.

In what follows, we denote by $B_R$ an open ball with radius $R$.  The symbol $\chi_X$ denotes the characteristic function
of a set $X$.

\subsection{Statement of the result}

Now let us describe the result more precisely.
First we find a solution of  (1.6) under special initial data described below.

Let $\eta>0$ and $R_2>R_1>0$.  Let $m=1,2,\cdots$.
In this paper, we choose two special initial data $\Phi(\,\cdot\,,0)=\Phi_{ext}(x), \Phi_{int}(x)$:
$$
\left\{
\begin{array}{l}
\displaystyle
\Phi_{ext,m}(x)
=
(\eta^2-\vert x-z\vert^2)^m\chi_{B_{\eta}}(x),\\
\\
\displaystyle
\Phi_{int,m}(x)
=(R_2^2-\vert x-z\vert^2)^m(\vert x-z\vert^2-R_1^2)^m\chi_{B_{R_2}\setminus B_{R_1}}(x),
\end{array}
\right.
$$
where ball $B_{\eta}$ and two concentric balls $B_{R_1}$ and $B_{R_2}$ are chosen in such a way that
$\overline{B_{\eta}}\cap\Omega=\emptyset$ and $\Omega\subset B_{R_1}$.  Here we made use of the common symbol $z$ to denote the center of $B_{\eta}$
in $\Phi_{ext,m}$ and of $B_{R_1}$ and $B_{R_2}$ in $\Phi_{int,m}$.
Note that both $\Phi_{ext,m}$ and $\Phi_{int,m}$ belong to $H^m(\Bbb R^3)\setminus H^{m+1}(\Bbb R^3)$.

Thus we have $\text{supp}\,\Phi(\,\cdot\,,0)=K$, where
$$\displaystyle
K
=
\left\{
\begin{array}{ll}
\displaystyle
\overline{B_{\eta}} & \text{if $\Phi(\,\cdot\,,0)=\Phi_{ext,m}$,}\\
\\
\displaystyle
\overline{B_{R_2}}\setminus B_{R_1} & \text{if $\Phi(\,\cdot\,,0)=\Phi_{int,m}$}
\end{array}
\right.
$$
and $K$ satisfies $K\cap\Omega=\emptyset$.

Using the Fourier transform method, we see that
there exists a unique $\Phi$ in the class
$C([0,\,T];H^m(\Bbb R^3))\cap H^1(0,\,T; H^{m-1}(\Bbb R^3))$
that satisfies (1.6) and the initial data $\Phi(\,\cdot\,,0)=\Phi_{\star,m}$, where $\star=\text{ext}$ or $\text{int}$.
Then choose $\mbox{\boldmath $v$}$ given by (1.7) which belongs to $H^1(0,\,T;H^{m-2}(\Bbb R^3))$
and thus $\mbox{\boldmath $f$}$ given by (1.4) belongs to $H^1(0,\,T;H^{m-2-\frac{1}{2}}(\partial\Omega))$.

We choose an arbitrary fixed $m\ge 4$.  
This implies $\mbox{\boldmath $f$}\in\,H^1(0,\,T;H^{\frac{3}{2}}(\partial\Omega))$.  
Then, it is known that there exists a unique solution $(\mbox{\boldmath $u$},p)$ of (1.1) such that
$\mbox{\boldmath $u$}\in L^2(0,\,T;H^2(\Omega\setminus\overline D)^3)
\cap C([0,\,T]; H^1(\Omega\setminus\overline D)^3)\cap H^1(0,\,T;L^2(\Omega\setminus\overline D)^3)$
and $p\in L^2(0,\,T;H^1(\Omega\setminus\overline D)/\Bbb R)$.
See Proposition 3.2 in \cite{MST}.

Now we state the main result of this paper.

\proclaim{\noindent Theorem 1.1.}

\noindent
(i)  Let $T$ be an arbitrary fixed positive number.
Then, there exists a positive number $\tau_0$ such that, for all $\tau\ge\tau_0$
$I_T(\tau;\mbox{\boldmath $v$},0)>0$ and we have
$$\displaystyle
\lim_{\tau\longrightarrow\infty}\frac{1}{\sqrt{\tau}}\log I_T(\tau;\mbox{\boldmath $v$},0)
= -2\sqrt{\frac{\rho}{\mu}}\text{dist}\,(D,K).
\tag {1.8}
$$

\noindent
(ii)  We have
$$\displaystyle
\lim_{\tau\longrightarrow\infty}e^{\sqrt{\tau} \,T}I_T(\tau;\mbox{\boldmath $v$},0)
=
\left\{\begin{array}{ll}
\displaystyle
\infty & 
\text{if $\displaystyle T>2\sqrt{\frac{\rho}{\mu}}\text{dist}\,(D,K)$,}
\\
\\
\displaystyle
0 & 
\text{if $\displaystyle T<2\sqrt{\frac{\rho}{\mu}}\text{dist}\,(D,K)$.}
\end{array}
\right.
$$

\noindent
(iii)  If $\displaystyle T=2\sqrt{\frac{\rho}{\mu}}\text{dist}\,(D,K)$, 
then we have, as $\tau\longrightarrow\infty$
$$\displaystyle
e^{\sqrt{\tau}\,T}\vert I_T(\tau;\mbox{\boldmath $v$},0)\vert=O(\tau^3).
\tag {1.9}
$$

\endproclaim

\begin{figure}[htbp]
\begin{center}
\epsfxsize=7cm
\epsfysize=5.5cm
\epsfbox{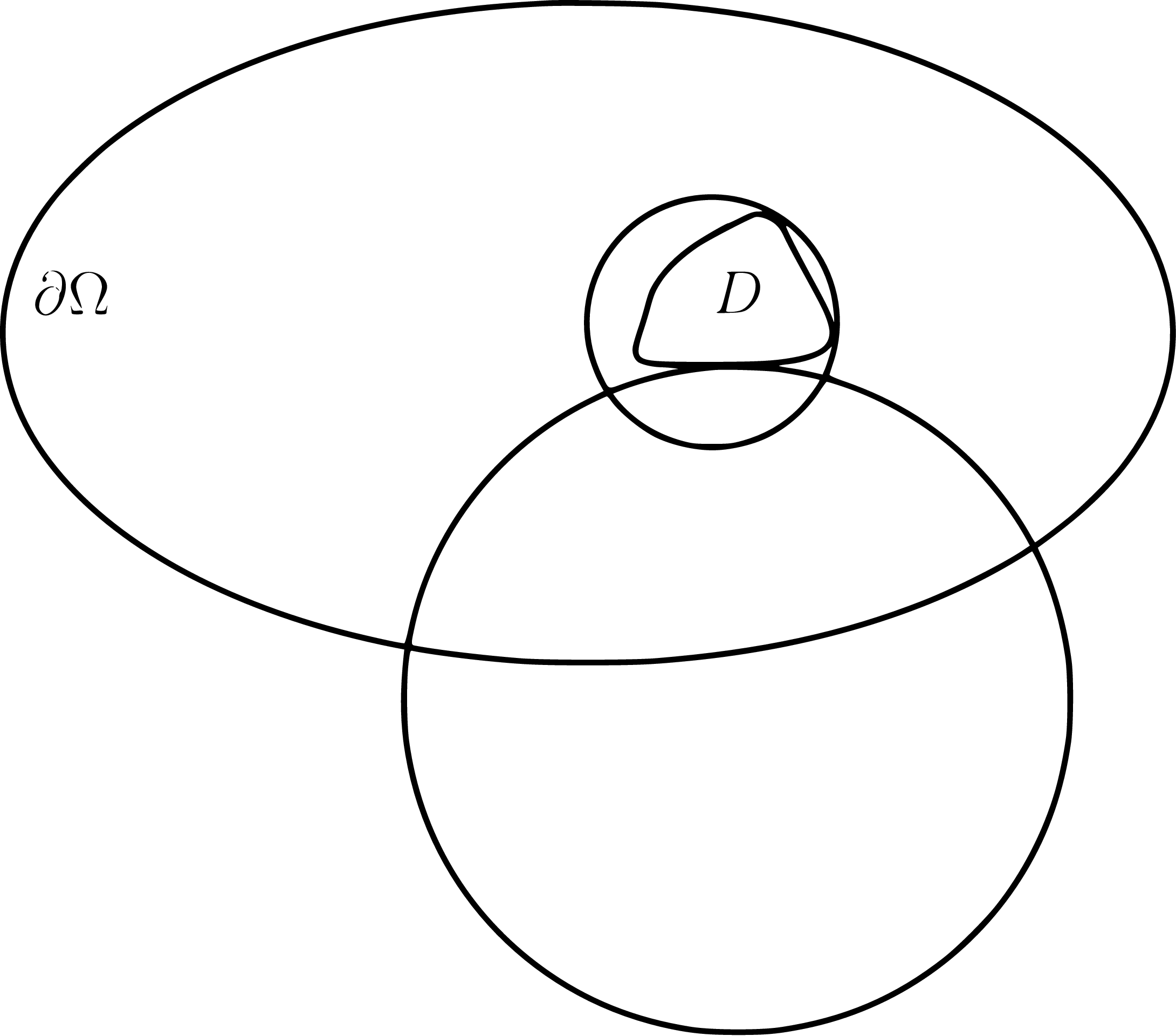}
\caption{An illustration of two large and small spheres centered at two different points
$z_1\in\Bbb R^3\setminus\overline\Omega$ and $z_2\in\Omega$ with radii $d_D(z_1)$ and $R_D(z_2)$, respectively.
}\label{fig1}
\end{center}
\end{figure}

It is easy to see that 
$$\displaystyle
\text{dist}\,(D,K)
=\left\{
\begin{array}{ll}
\displaystyle
d_D(z_1)-\eta, & \text{if $K=\overline{B_{\eta}}$,}\\
\\
\displaystyle
R_1-R_D(z_2), & \text{if $K=\overline{B_{R_2}}\setminus B_{R_1}$,}
\end{array}
\right.
$$
where $z_1$ denotes the center point of $B_{\eta}$ and $z_2$ the common center point of $B_{R_1}$ and $B_{R_2}$; 
$$\begin{array}{ll}
\displaystyle
d_D(z_1)=\text{dist}\,(\{z_1\},D), & R_D(z_2)=\sup_{x\in D}\vert x-z_2\vert.
\end{array}
$$
Thus by Theorem 1.1 one can extract two quantities $d_D(z_1)$ and $R_D(z_2)$
which yield the largest/smallest sphere whose exterior/interior contains $D$
from the indicator function according to the two choices of $\Phi(\,\cdot\,,0)$.
See Figure \ref{fig1} for an illustration of two spheres
centered at two different points
$z_1$ and $z_2$ with radii $d_D(z_1)$ and $R_D(z_2)$, respectively.

\vspace{0.0cm}

It should be emphasized that there is no restriction on $T$.  This indicates that the propagation
speed of the signal governed by the Stokes system is not finite similar to that of the signal
governed by the heat equation.

$\quad$

{\bf\noindent Remark 1.1.}
Since $T$ in (i) is {\it arbitrary}, (ii) can be easily deduced from (i).  Thus the point is to establish
(i) and (iii).

$\quad$

{\bf\noindent Remark 1.2.}  Mathematically, from (1.8) we obtain $\text{dist}\,(D,K)$ from the data over an arbitrary fixed time
interval $]0,\,T[$.  So this shows that the propagation speed of the signal governed by the Stokes system is infinite.

$\quad$

This paper is concerned with the reconstruction issue for an obstacle embedded in a {\it bounded domain}.  
For the uniqueness and stability issue see \cite{ACFKH} including the stationary case.
In \cite{B} a stability estimate of log-log type in the stationary case has been established.
For a numerical reconstruction method in the stationary case based on the method of fundamental solutions see
\cite{AKS}.  In \cite{BD} a combination of the quasi-reversibility method and the level set one has been applied.

In the article \cite{HUW}, the idea of the original enclosure method \cite{E00, E01} has been applied to the Stokes system in the stationary case (with {\it variable viscosity}).  They state an extraction procedure of the distance of an arbitrary point $x_0$ {\it outside the convex hull}
of $\overline\Omega$ to an unknown obstacle $D$ from the associated Dirichlet-to-Neumann map acting on
infinitely many special Dirichlet data.  
Their Dirichlet data are given by the trace onto $\partial\Omega$ of  so-called complex spherical wave 
solutions of the stationary Stokes system
in an open neighbourhood of $\overline\Omega$ having several parameters.

Besides, we mention an article using the idea of the enclosure method in {\it the time domain}.
In \cite{MST} the response operator acting infinitely many $\mbox{\boldmath $f$}$ has been used.
They are the separation of
variables type having the form
$$\begin{array}{ll}
\displaystyle
\mbox{\boldmath $f$}(x,t)=\chi(t)\mbox{\boldmath $v$}(x), & (x,t)\in\partial\Omega\times\,]0,\,T[,
\end{array}
\tag {1.10}
$$
where the function $\mbox{\boldmath $v$}$ defined in an open neighbourhood $\tilde{\Omega}$ of $\overline\Omega$ 
and satisfies the following system depending on a large parameter $\tau$ (in our notation) with an appropriate function $q$:
$$\left\{
\begin{array}{ll}
\displaystyle
\text{div}\,\sigma(\mbox{\boldmath $v$},q)-\tau\rho\mbox{\boldmath $v$}=\mbox{\boldmath $0$}, & x\in\tilde{\Omega}
\\
\\
\displaystyle
\nabla\cdot\mbox{\boldmath $v$}=0, & x\in\tilde{\Omega}.
\end{array}
\right.
\tag {1.11}
$$
The function $\chi(t)$ also depends $\tau>0$ and satisfies
$\chi(0)=0$ and $\chi(t)>0$ for all $0<t\le T$ with the normalized condition $\int_0^Te^{-\tau t}\,\chi(t)\,dt=1$.
This is just under the framework of an earlier version of the time domain enclosure method with infinitely many input proposed in \cite{IFR}.

In \cite{MST} they used two types of the solutions.  One is the real exponential solution having the form
$$\begin{array}{ll}
\displaystyle
\mbox{\boldmath $v$}(x)=e^{\tilde{\tau}\,\mbox{\boldmath $m$}\cdot x}\,\mbox{\boldmath $l$}, & x\in\Bbb R^3
\end{array}
\tag{1.12}
$$
with $q=0$ and $\tilde{\tau}=\sqrt{\tau}\sqrt{\rho/\mu}$, where $\mbox{\boldmath $l$}$ and $\mbox{\boldmath $m$}$ 
are two unit constant vectors satisfying $\mbox{\boldmath $l$}\cdot\mbox{\boldmath $m$}=0$.

Another one is a {\it singular solution} having its singularity on a line not intersecting with the convex hull of
$\overline\Omega$. 
The solution has the following form.  Let $x_0$ be an arbitrary point
outside of the convex hull of $\overline\Omega$. 
Choose a unit vector $\mbox{\boldmath $c$}$ in such a way that the line 
$l:x=x_0+\lambda\,\mbox{\boldmath $c$}$,
$-\infty<\lambda<\infty$ lies in the outside of the convex hull of $\overline\Omega$.
Then, the $\mbox{\boldmath $v$}$ takes the form
$$\begin{array}{ll}
\displaystyle
\mbox{\boldmath $v$}(x)=e^{-\tilde{\tau}\vert x-x_0\vert}
\frac{\mbox{\boldmath $c$}\times(x-x_0)}
{\vert\mbox{\boldmath $c$}\times(x-x_0)\vert^2},
&
x\in\Bbb R^3\setminus l
\end{array}
\tag {1.13}
$$
and $q=0$.

Using a spherical coordinates depending on $\mbox{\boldmath $c$}$ of $\Bbb R^3$, 
they showed that the $\mbox{\boldmath $v$}$ together with $q=0$ satisfies (1.11)
and essentially it's energy integral locally outside the line $l$ behaves like that of the function
$$\displaystyle
e^{-\tilde{\tau}\,\vert x-x_0\vert}.
$$
This means that the essential part is reduced to the heat equation case considered in \cite{IK1}.

Using $\mbox{\boldmath $v$}$ given by (1.12) and (1.13) 
with $q=0$, they prescribe $\mbox{\boldmath $f$}$ given by (1.10)
and gave extraction formulae of the quantities $h_D(\mbox{\boldmath $m$})=\sup_{x\in\,D}x\cdot\mbox{\boldmath $m$}$ and  
$\text{dist}(x_0,D)$ provided $x_0$ is placed outside of the convex hull
of $\overline\Omega$.  In our method $\mbox{\boldmath $f$}$ is independent of $\tau$ and yields
$\text{dist}(x_0,D)$ for all $x_0\in\Bbb R^3\setminus\overline\Omega$ and the quantity
$R_D(x_0)$ for all $x_0\in\Bbb R^3$ which is not extracted in their paper.

\subsection{Outline of the paper}

A brief outline of this paper is as follows. Theorems 1.1
is proved in Section 3. The proof is based on:
a representation formula of the indicator function together with some energy estimates
which are proved in Section 2;
Proposition 3.1 which gives explicit 
formulae in some domains for the solutions of the modified Helmholtz equation with special inhomogeneous terms
in the whole space and the proof is given in Section 4.

\section{Preliminaries}

In this section, the symbol $\mbox{\boldmath $n$}$ denotes also the unit outward normal to $\partial D$.

The function $\mbox{\boldmath $v$}$ is given by a solution of (1.3) with
$\text{supp}\,\mbox{\boldmath $v$}(\,\cdot\,,0)\cap\Omega=\emptyset$
and $\mbox{\boldmath $f$}$ given by (1.4).  Note that the form of $\mbox{\boldmath $v$}$
is not specified unlike (1.7).

\subsection{Representation formula and its direct consequence}

From (1.1) we have
$$\left\{
\begin{array}{ll}
\displaystyle
\text{div}\,\sigma(\mbox{\boldmath $w$},\tilde{p})-\tau\rho\mbox{\boldmath $w$}
=e^{-\tau T}\mbox{\boldmath $F$},
&
\displaystyle
x\in\Omega\setminus\overline D,
\\
\\
\displaystyle
\nabla\cdot\mbox{\boldmath $w$}=0,
&
\displaystyle
x\in\Omega\setminus\overline D,
\\
\\
\displaystyle
\mbox{\boldmath $w$}=\mbox{\boldmath $0$},
&
\displaystyle
x\in\partial D,\\
\\
\displaystyle
\mbox{\boldmath $w$}=\mbox{\boldmath $w$}_0,
&
\displaystyle
x\in\partial\Omega,
\end{array}
\right.
\tag {2.1}
$$
where
$$\displaystyle
\mbox{\boldmath $F$}(x)=\rho \mbox{\boldmath $u$}(x,T).
$$
And also from (1.3) we have
$$\left\{
\begin{array}{ll}
\displaystyle
\text{div}\,\sigma(\mbox{\boldmath $w$}_0,\tilde{q})-\tau\rho\mbox{\boldmath $w$}_0
+\mbox{\boldmath $v$}(x,0)
=e^{-\tau T}\mbox{\boldmath $F$}_0,
&
\displaystyle
x\in\Bbb R^3,
\\
\\
\displaystyle
\nabla\cdot\mbox{\boldmath $w$}_0=0,
&
\displaystyle
x\in\Bbb R^3,
\end{array}
\right.
\tag {2.2}
$$
where
$$\displaystyle
\mbox{\boldmath $F$}_0(x)=\rho \mbox{\boldmath $v$}(x,T).
$$
Noting that $\text{supp}\,\mbox{\boldmath $v$}(\,\cdot\,,0)\cap\Omega=\emptyset$, from (2.1) and (2.2) we obtain
$$\begin{array}{l}
\displaystyle
\,\,\,\,\,\,
\int_{\partial\Omega}
\left(\sigma(\mbox{\boldmath $w$},\tilde{p})\mbox{\boldmath $n$}\cdot\mbox{\boldmath $w$}_0
-\sigma(\mbox{\boldmath $w$}_0,\tilde{q})\mbox{\boldmath $n$}\cdot\mbox{\boldmath $w$}\right)\,dS\\
\\
\displaystyle
=\int_{\partial D}\sigma(\mbox{\boldmath $w$},\tilde{p})\mbox{\boldmath $n$}\cdot\mbox{\boldmath $w$}_0\,dS\\
\\
\displaystyle
\,\,\,
+
\int_{\Omega\setminus\overline D}
\left(\sigma(\mbox{\boldmath $w$},\tilde{p})\cdot\nabla\mbox{\boldmath $w$}_0
-\sigma(\mbox{\boldmath $w$}_0,\tilde{q})\cdot\nabla\mbox{\boldmath $w$}\right)\,dx\\
\\
\displaystyle
\,\,\,
+e^{-\tau T}
\int_{\Omega\setminus\overline D}
(\mbox{\boldmath $F$}\cdot\mbox{\boldmath $w$}_0
-\mbox{\boldmath $F$}_0\cdot\mbox{\boldmath $w$})\,dx.
\end{array}
\tag {2.3}
$$
Using $\nabla\cdot\mbox{\boldmath $w$}_0=\nabla\cdot\mbox{\boldmath $w$}=\mbox{\boldmath $0$}$ in $\Omega\setminus\overline D$
and
$$
\displaystyle
\text{Sym}\,\nabla\mbox{\boldmath $w$}_0\cdot\nabla\mbox{\boldmath $w$}
=\text{Sym}\,\nabla\mbox{\boldmath $w$}\cdot\nabla\mbox{\boldmath $w$}_0,
$$
from (2.3) we obtain the first expression
of the indicator function:
$$\begin{array}{l}
\displaystyle
\,\,\,\,\,\,
I_T(\tau;\mbox{\boldmath $v$},q)\\
\\
\displaystyle
=\int_{\partial D}\sigma(\mbox{\boldmath $w$},\tilde{p})\mbox{\boldmath $n$}\cdot\mbox{\boldmath $w$}_0\,dS
+e^{-\tau T}
\int_{\Omega\setminus\overline D}
(\mbox{\boldmath $F$}\cdot\mbox{\boldmath $w$}_0
-\mbox{\boldmath $F$}_0\cdot\mbox{\boldmath $w$})\,dx.
\end{array}
\tag {2.4}
$$

Set
$$\displaystyle
\mbox{\boldmath $R$}=\mbox{\boldmath $w$}-\mbox{\boldmath $w$}_0.
$$

We further rewrite formula (2.4).

\proclaim{\noindent Proposition 2.1.}
We have
$$\displaystyle
I_T(\tau;\mbox{\boldmath $v$},q)=J(\tau)+E(\tau)+e^{-\tau T}{\cal R}(\tau),
\tag {2.5}
$$
where
$$\displaystyle
J(\tau)
=\int_{D}(2\mu\,\vert\text{Sym}\,\nabla\mbox{\boldmath $w$}_0\vert^2+\tau\rho\vert\mbox{\boldmath $w$}_0\vert^2)\,dx,
\tag {2.6}
$$
$$\displaystyle
E(\tau)
=\int_{\Omega\setminus\overline D}(2\mu\,\vert\text{Sym}\,\nabla\mbox{\boldmath $R$}\vert^2+\tau\rho\vert\mbox{\boldmath $R$}\vert^2)\,dx
\tag {2.7}
$$
and
$$\displaystyle
{\cal R}(\tau)=\int_D\mbox{\boldmath $F$}_0\cdot\mbox{\boldmath $w$}_0\,dx
+\int_{\Omega\setminus\overline D}(\mbox{\boldmath $F$}-2\mbox{\boldmath $F$}_0)\cdot\mbox{\boldmath $R$}\,dx
+\int_{\Omega\setminus\overline D}
(\mbox{\boldmath $F$}-
\mbox{\boldmath $F$}_0)\cdot\mbox{\boldmath $w$}_0\,dx.
\tag {2.8}
$$

\endproclaim

{\it\noindent Proof.}
From (2.1) and (2.2) we have
$$\left\{
\begin{array}{ll}
\displaystyle
\text{div}\,\sigma(\mbox{\boldmath $R$},r)-\tau\rho\mbox{\boldmath $R$}
=e^{-\tau T}(\mbox{\boldmath $F$}-\mbox{\boldmath $F$}_0),
&
\displaystyle
x\in\Omega\setminus\overline D,
\\
\\
\displaystyle
\nabla\cdot\mbox{\boldmath $R$}=0,
&
\displaystyle
x\in\Omega\setminus\overline D,
\\
\\
\displaystyle
\mbox{\boldmath $R$}=-\mbox{\boldmath $w$}_0,
&
\displaystyle
x\in\partial D,\\
\\
\displaystyle
\mbox{\boldmath $R$}=\mbox{\boldmath $0$},
&
\displaystyle
x\in\partial\Omega,
\end{array}
\right.
\tag {2.9}
$$
where
$$
\displaystyle
r=\tilde{p}-\tilde{q}.
$$

Write
$$\begin{array}{l}
\displaystyle
\,\,\,\,\,\,
\int_{\partial D}\sigma(\mbox{\boldmath $w$},\tilde{p})\mbox{\boldmath $n$}\cdot\mbox{\boldmath $w$}_0\,dS
\\
\\
\displaystyle
=\int_{\partial D}\sigma(\mbox{\boldmath $w$}_0,\tilde{q})\mbox{\boldmath $n$}\cdot\mbox{\boldmath $w$}_0\,dS
+\int_{\partial D}\sigma(\mbox{\boldmath $R$}, r)\mbox{\boldmath $n$}\cdot\mbox{\boldmath $w$}_0\,dS
\end{array}
$$
It follows from (2.2) in $D$ that
$$
\displaystyle
\int_{\partial D}\sigma(\mbox{\boldmath $w$}_0,\tilde{q})\mbox{\boldmath $n$}\cdot\mbox{\boldmath $w$}_0\,dS
=\int_D\left(\sigma(\mbox{\boldmath $w$}_0,\tilde{q})\cdot\nabla\mbox{\boldmath $w$}_0+\tau\rho\vert\mbox{\boldmath $w$}_0\vert^2\right)\,dx
+e^{-\tau T}\int_D\mbox{\boldmath $F$}_0\cdot\mbox{\boldmath $w$}_0\,dx.
$$
It follows from (2.9) that
$$\begin{array}{l}
\displaystyle
\,\,\,\,\,\,
\int_{\partial D}\sigma(\mbox{\boldmath $R$},r)\mbox{\boldmath $n$}\cdot\mbox{\boldmath $w$}_0\,dS
\\
\\
\displaystyle
=-\int_{\partial D}\sigma(\mbox{\boldmath $R$}, r)\mbox{\boldmath $n$}\cdot\mbox{\boldmath $R$}\,dS\\
\\
\displaystyle
=\int_{\partial\,(\Omega\setminus\overline D)}\sigma(\mbox{\boldmath $R$}, r)\mbox{\boldmath $n$}\cdot\mbox{\boldmath $R$}\,dS
\\
\\
\displaystyle
=\int_{\Omega\setminus\overline D}(\sigma(\mbox{\boldmath $R$}, r)\cdot\nabla\mbox{\boldmath $R$}+\tau\rho\vert\mbox{\boldmath $R$}\vert^2)\,dx
+e^{-\tau T}\int_{\Omega\setminus\overline D}(\mbox{\boldmath $F$}-\mbox{\boldmath $F$}_0)\cdot\mbox{\boldmath $R$}dx.
\end{array}
\tag {2.10}
$$
Thus we obtain
$$\begin{array}{l}
\displaystyle
\,\,\,\,\,\,
\int_{\partial D}\sigma(\mbox{\boldmath $w$},\tilde{p})\mbox{\boldmath $n$}\cdot\mbox{\boldmath $w$}_0\,dS
\\
\\
\displaystyle
=\int_{D}(\sigma(\mbox{\boldmath $w$}_0,\tilde{q})\cdot\nabla\mbox{\boldmath $w$}_0+\tau\rho\vert\mbox{\boldmath $w$}_0\vert^2)\,dx
+\int_{\Omega\setminus\overline D}(\sigma(\mbox{\boldmath $R$},r)\cdot\nabla\mbox{\boldmath $R$}+\tau\rho\vert\mbox{\boldmath $R$}\vert^2)\,dx
\\
\\
\displaystyle
\,\,\,
+e^{-\tau T}\left\{
\int_D\mbox{\boldmath $F$}_0\cdot\mbox{\boldmath $w$}_0\,dx+\int_{\Omega\setminus\overline D}
(\mbox{\boldmath $F$}-
\mbox{\boldmath $F$}_0)\cdot\mbox{\boldmath $R$}dx
\right\}.
\end{array}
\tag {2.11}
$$
Since we have $\nabla\cdot\mbox{\boldmath $w$}_0=0$ in $D$ and $\nabla\cdot\mbox{\boldmath $w$}=0$ in $\Omega\setminus\overline D$,
using the identity $\text{Sym}\,A\cdot A=\vert\text{Sym}\,A\vert^2$, we obtain the expression
$$\left\{\begin{array}{l}
\displaystyle
\int_{D}(\sigma(\mbox{\boldmath $w$}_0,\tilde{q})\cdot\nabla\mbox{\boldmath $w$}_0+\tau\rho\vert\mbox{\boldmath $w$}_0\vert^2)\,dx=J(\tau)\\
\\
\displaystyle
\int_{\Omega\setminus\overline D}(\sigma(\mbox{\boldmath $R$},r)\cdot\nabla\mbox{\boldmath $R$}+\tau\rho\vert\mbox{\boldmath $R$}\vert^2)\,dx=E(\tau).
\end{array}
\right.
\tag {2.12}
$$
Moreover, one has
$$\begin{array}{l}
\,\,\,\,\,\,
\displaystyle
\int_{\Omega\setminus\overline D}
(\mbox{\boldmath $F$}-
\mbox{\boldmath $F$}_0)\cdot\mbox{\boldmath $R$}dx
+\int_{\Omega\setminus\overline D}
(\mbox{\boldmath $F$}\cdot\mbox{\boldmath $w$}_0
-\mbox{\boldmath $F$}_0\cdot\mbox{\boldmath $w$})\,dx
\\
\\
\displaystyle
=
\int_{\Omega\setminus\overline D}(\mbox{\boldmath $F$}-2\mbox{\boldmath $F$}_0)\cdot\mbox{\boldmath $R$}\,dx
+\int_{\Omega\setminus\overline D}
(\mbox{\boldmath $F$}-
\mbox{\boldmath $F$}_0)\cdot\mbox{\boldmath $w$}_0\,dx.
\end{array}
$$
Substituting (2.11) into (2.4) and using this together with (2.12), we obtain (2.5).

\noindent
$\Box$

Now we state two estimates which tell us that the asymptotic behaviour 
of the terms  $E(\tau)$ and $e^{-\tau T}{\cal R}(\tau)$ as $\tau\rightarrow\infty$
can be controlled from above by that of $J(\tau)$ and $\Vert\mbox{\boldmath $w$}_0\Vert_{L^2(\Omega)}$.

\proclaim{\noindent Proposition 2.2.}
We have, as $\tau\longrightarrow\infty$
$$\displaystyle
E(\tau)=O(\tau J(\tau)+e^{-2\tau T})
\tag {2.13}
$$
and
$$\displaystyle
e^{-\tau T}\vert{\cal R}(\tau)\vert
=O(e^{-\tau T}\Vert \mbox{\boldmath $w$}_0\Vert_{L^2(\Omega)}
+e^{-\tau T}J(\tau)^{1/2}+\tau^{-1/2}e^{-2\tau T}).
\tag {2.14}
$$

\endproclaim

{\it\noindent Proof.} 
First we give a proof of (2.13).
From (2.10) and the second formula on (2.12),
we have
$$\begin{array}{l}
\displaystyle
\,\,\,\,\,\,
\int_{\Omega\setminus\overline D}
\left(2\mu\vert\text{Sym}\,\nabla\mbox{\boldmath $R$}\vert^2
+\tau\rho\left\vert\mbox{\boldmath $R$}+e^{-\tau T}\frac{\mbox{\boldmath $F$}-\mbox{\boldmath $F$}_0}{2\tau\rho}\right\vert^2\right)\,dx
\\
\\
\displaystyle
=\int_{\partial D}\sigma(\mbox{\boldmath $R$},r)\mbox{\boldmath $n$}\cdot\mbox{\boldmath $w$}_0\,dS
+
\frac{e^{-2\tau T}}{4\tau\rho}
\int_{\Omega\setminus\overline D}\vert\mbox{\boldmath $F$}-\mbox{\boldmath $F$}_0\vert^2\,dx.
\end{array}
$$
This together with (2.7) immediately yields
$$\displaystyle
E(\tau)
\le
2\int_{\partial D}\sigma(\mbox{\boldmath $R$},r)\mbox{\boldmath $n$}\cdot\mbox{\boldmath $w$}_0\,dS
+O(\tau^{-1}e^{-2\tau T}).
\tag {2.15}
$$
Here we choose $\tilde{\mbox{\boldmath $w$}_0}\in H^1(\Omega\setminus\overline D)^3$ in such a way that
$$\left\{
\begin{array}{ll}
\displaystyle
\text{div}\,\tilde{\mbox{\boldmath $w$}_0}=0, & x\in\Omega\setminus\overline D,\\
\\
\displaystyle
\tilde{\mbox{\boldmath $w$}_0}=\mbox{\boldmath $0$}, & x\in\partial\Omega,\\
\\
\displaystyle
\tilde{\mbox{\boldmath $w$}_0}=\mbox{\boldmath $w$}_0, & x\in\partial D
\end{array}
\right.
$$
and
$$\displaystyle
\Vert\tilde{\mbox{\boldmath $w$}_0}\Vert_{H^1(\Omega\setminus\overline D)}
\le C
\Vert\mbox{\boldmath $w$}_0\Vert_{H^{1/2}(\partial D)}
\le C'
\Vert\mbox{\boldmath $w$}_0\Vert_{H^1(D)}.
\tag {2.16}
$$
See the proof of Lemma 3.1 in \cite{MST} for this choice which is taken from \cite{Ga}.

Then, we can write
$$\begin{array}{l}
\,\,\,\,\,\,
\displaystyle
\int_{\partial D}\sigma(\mbox{\boldmath $R$},r)\mbox{\boldmath $n$}\cdot\mbox{\boldmath $w$}_0\,dS
\\
\\
\displaystyle
=\int_{\partial D}\sigma(\mbox{\boldmath $R$},r)\mbox{\boldmath $n$}\cdot\tilde{\mbox{\boldmath $w$}_0}\,dS\\
\\
\displaystyle
=-\int_{\partial(\Omega\setminus\overline D)}\sigma(\mbox{\boldmath $R$},r)\mbox{\boldmath $n$}\cdot\tilde{\mbox{\boldmath $w$}_0}\,dS.
\end{array}
\tag {2.17}
$$
Integration by parts and the first equation on (2.9) give
$$\begin{array}{l}
\,\,\,\,\,\,
\displaystyle
\int_{\partial(\Omega\setminus\overline D)}\sigma(\mbox{\boldmath $R$},r)\mbox{\boldmath $n$}\cdot\tilde{\mbox{\boldmath $w$}_0}\,dS
\\
\\
\displaystyle
=\int_{\Omega\setminus\overline D}\text{div}\,\sigma(\mbox{\boldmath $R$},r)\cdot\tilde{\mbox{\boldmath $w$}_0}\,dx
+\int_{\Omega\setminus\overline D}\sigma(\mbox{\boldmath $R$},r)\cdot\nabla\tilde{\mbox{\boldmath $w$}_0}\,dx\\
\\
\displaystyle
=\int_{\Omega\setminus\overline D}\text{div}\,\sigma(\mbox{\boldmath $R$},r)\cdot\tilde{\mbox{\boldmath $w$}_0}\,dx
+\int_{\Omega\setminus\overline D}2\mu\,\text{Sym}\,\nabla\mbox{\boldmath $R$}\cdot\text{Sym}\,\nabla\tilde{\mbox{\boldmath $w$}_0}\,dx\\
\\
\displaystyle
=\int_{\Omega\setminus\overline D}\tau\rho\mbox{\boldmath $R$}\cdot\tilde{\mbox{\boldmath $w$}_0}\,dx
+e^{-\tau T}\int_{\Omega\setminus\overline D}(\mbox{\boldmath $F$}-\mbox{\boldmath $F$}_0)\cdot\tilde{\mbox{\boldmath $w$}_0}\,dx
+\int_{\Omega\setminus\overline D}2\mu\,\text{Sym}\,\nabla\mbox{\boldmath $R$}\cdot\text{Sym}\,\nabla\tilde{\mbox{\boldmath $w$}_0}\,dx.
\end{array}
$$
From (2.16) we see that this right-hand side has the bound
$$\displaystyle
C\Vert \mbox{\boldmath $w$}_0\Vert_{H^1(D)}(\tau\Vert\mbox{\boldmath $R$}\Vert_{L^2(\Omega\setminus\overline D)}
+\Vert\text{Sym}\,\nabla\mbox{\boldmath $R$}\Vert_{L^2(\Omega\setminus\overline D)}+e^{-\tau T}).
$$
Applying this and (2.17) to the first term on (2.15), we obtain
$$\displaystyle
E(\tau)
\le 
C\Vert \mbox{\boldmath $w$}_0\Vert_{H^1(D)}(\tau\Vert\mbox{\boldmath $R$}\Vert_{L^2(\Omega\setminus\overline D)}
+\Vert\text{Sym}\,\nabla\mbox{\boldmath $R$}\Vert_{L^2(\Omega\setminus\overline D)}+e^{-\tau T})
+O(\tau^{-1}e^{-2\tau T}).
$$
Moreover, from (2.7) we have, for all $\tau\ge 1$
$$\displaystyle
\tau\Vert\mbox{\boldmath $R$}\Vert_{L^2(\Omega\setminus\overline D)}
+\Vert\text{Sym}\,\nabla\mbox{\boldmath $R$}\Vert_{L^2(\Omega\setminus\overline D)}
\le C \tau^{1/2} E(\tau)^{1/2}.
\tag {2.18}
$$
From these we obtain
$$\displaystyle
E(\tau)
\le 
C\Vert \mbox{\boldmath $w$}_0\Vert_{H^1(D)}
(\tau^{1/2} E(\tau)^{1/2}+e^{-\tau T})
+O(\tau^{-1}e^{-2\tau T}).
\tag {2.19}
$$
Here recall Korn's second inequality \cite{DuL}
$$\displaystyle
\Vert\mbox{\boldmath $w$}_0\Vert_{H^1(D)}
\le C(\Vert\mbox{\boldmath $w$}_0\Vert_{L^2(D)}^2+
\Vert\text{Sym}\,\nabla\mbox{\boldmath $w$}_0\Vert_{L^2(D)}^2)^{1/2}.
$$
Since this right-hand side has the bound $J(\tau)^{1/2}$, from (2.19)
we obtain
$$\displaystyle
E(\tau)
\le 
CJ(\tau)^{1/2}
(\tau^{1/2} E(\tau)^{1/2}+e^{-\tau T})
+O(\tau^{-1}e^{-2\tau T}).
$$
Then, a standard argument yields (2.13).

Next we give a proof of (2.14).
From (2.8) we have
$$\displaystyle
e^{-\tau T}\vert{\cal R}(\tau)\vert\le C(e^{-\tau T}\Vert\mbox{\boldmath $w$}_0\Vert_{L^2(\Omega)}+
e^{-\tau T}\Vert\mbox{\boldmath $R$}\Vert_{L^2(\Omega\setminus\overline D)}).
$$
Besides (2.18) yields the estimate
$$
\displaystyle
\Vert\mbox{\boldmath $R$}\Vert_{L^2(\Omega\setminus\overline D)}
\le C\tau^{-1/2}E(\tau)^{1/2}.
$$
Now from these together with (2.13) we obtain (2.14).

\noindent
$\Box$

\section{Proof of Theorem 1.1}

In this section the form of $\mbox{\boldmath $v$}$
is specified as (1.7).

\subsection{Propagation estimates}

The following two lemmas are concernd with the asymptotic behaviour of 
$J(\tau)$ and ${\cal R}(\tau)$ as $\tau\longrightarrow\infty$.

\proclaim{\noindent Lemma 3.1.}
Let $U$ be an arbitrarly bounded open subset of $\Bbb R^3$.
We have, as $\tau\longrightarrow\infty$
$$
\displaystyle
\sqrt{\tau}\Vert\mbox{\boldmath $w$}_0\Vert_{L^2(U)}
+\Vert\text{Sym}\,\nabla\mbox{\boldmath $w$}_0\Vert_{L^2(U)}
=O(\tau e^{-\sqrt{\tau}\,\sqrt{\rho/\mu}\,\text{dist}\,(U,K)}+\tau^{-1/2}e^{-\tau T})
\tag {3.1}
$$
\endproclaim
{\it\noindent Proof.}
Let $\mbox{\boldmath $v$}_{00}\in H^1(\Bbb R^3)^3$ be the unique weak solution of 
$$\begin{array}{ll}
\displaystyle
\mu\Delta\mbox{\boldmath $v$}-\tau\rho\mbox{\boldmath $v$}+\Phi(x,0)\mbox{\boldmath $a$}
=\mbox{\boldmath $0$}, &
x\in\Bbb R^3.
\end{array}
\tag {3.2}
$$
We have the expression
$$\displaystyle
\mbox{\boldmath $v$}_{00}
(x)
=\left(\frac{1}{4\pi\mu}\int_{\Bbb R^3}
\frac{e^{-\sqrt{\tau}\sqrt{\rho/\mu}\,\vert x-y\vert}}{\vert x-y\vert}\,\Phi(\,\cdot\,,0)\,dy\right)\,\mbox{\boldmath $a$}.
\tag {3.3}
$$
Since $\Phi(\,\cdot\,,0)\in H^m(\Bbb R^3)$ we have $\mbox{\boldmath $v$}_{00}\in H^{m+2}(\Bbb R^3)^2$.

Define
$$\displaystyle
\mbox{\boldmath $w$}_{00}=\nabla\times\mbox{\boldmath $v$}_{00}\in H^{m+1}(\Bbb R^3)^3.
\tag {3.4}
$$
Then the pair $(\mbox{\boldmath $w$}_{00},\tilde{q})$ with $\tilde{q}=0$ satisfies
$$\left\{
\begin{array}{ll}
\displaystyle
\text{div}\,\sigma(\mbox{\boldmath $w$}_{00},\tilde{q})-\tau\rho\mbox{\boldmath $w$}_{00}
+\nabla\times(\Phi(x,0)\mbox{\boldmath $a$})
=\mbox{\boldmath $0$},
&
\displaystyle
x\in\Bbb R^3,
\\
\\
\displaystyle
\nabla\cdot\mbox{\boldmath $w$}_{00}=0,
&
\displaystyle
x\in\Bbb R^3.
\end{array}
\right.
\tag {3.5}
$$

From the expression (3.3) and (3.4) we have
$$\displaystyle
\sqrt{\tau}\,\Vert\mbox{\boldmath $w$}_{00}\Vert_{L^2(U)}
+\Vert\nabla\mbox{\boldmath $w$}_{00}\Vert_{L^2(U)}
\le O(\tau e^{-\sqrt{\tau}\,\sqrt{\rho/\mu}\,\text{dist}\,(U,K)}).
\tag {3.6}
$$
Set
$$\displaystyle
\mbox{\boldmath $\epsilon$}_0=e^{-\tau T}(\mbox{\boldmath $w$}_{0}-\mbox{\boldmath $w$}_{00}).
$$
We have
$$\displaystyle
\mbox{\boldmath $w$}_0=\mbox{\boldmath $w$}_{00}+e^{-\tau T}\mbox{\boldmath $\epsilon$}_0
$$
and it follows from (2.2), (3.2) and (3.5) that $\mbox{\boldmath $\epsilon$}_0$ satisfies
$$\left\{
\begin{array}{ll}
\displaystyle
\text{div}\,\sigma(\mbox{\boldmath $\epsilon$}_0,\tilde{q})-\tau\rho\mbox{\boldmath $\epsilon$}_0
=\mbox{\boldmath $F$}_0,
&
\displaystyle
x\in\Bbb R^3,\\
\\
\displaystyle
\nabla\cdot\mbox{\boldmath $\epsilon$}_0=0, & x\in\Bbb R^3.
\end{array}
\right.
$$
Since we have
$$\displaystyle
\int_{\Bbb R^3}(2\mu\vert\text{Sym}\,\nabla\mbox{\boldmath $\epsilon$}_0\vert^2
+\tau\rho\vert\mbox{\boldmath $\epsilon$}_0\vert^2+\mbox{\boldmath $F$}_0\cdot\mbox{\boldmath $\epsilon$}_0)\,dx=0,
$$
one gets
$$\displaystyle
\int_{\Bbb R^3}(2\mu\vert\text{Sym}\,\nabla\mbox{\boldmath $\epsilon$}_0\vert^2+\tau\rho\vert\mbox{\boldmath $\epsilon$}_0\vert^2)\,dx
\le 2\times 2\times\frac{\Vert\mbox{\boldmath $F$}_0\Vert_{L^2(\Bbb R^3)}^2}{4\tau\rho}=O(\tau^{-1}).
$$
This gives
$$\displaystyle
\sqrt{\tau}\Vert\mbox{\boldmath $\epsilon$}_0\Vert_{L^2(\Bbb R^3)}
+\Vert\text{Sym}\,\nabla\mbox{\boldmath $\epsilon$}_0\Vert_{L^2(\Bbb R^3)}
=O(\tau^{-1/2}).
\tag {3.7}
$$
Then, from (3.6) and (3.7) we obtain (3.1).

\noindent
$\Box$

To obtain a lower estimate for $J(\tau)$ we need a detailed expression of $\mbox{\boldmath $w$}_{00}=\nabla\times\mbox{\boldmath $v$}_{00}$
in a neighbouhood of $\overline D$.

Set
$$
\displaystyle
\mbox{\boldmath $v$}_{00}=
\left\{
\begin{array}{ll}
\displaystyle
\mbox{\boldmath $v$}_{ext,m}, & \text{if $\Phi(\,\cdot\,,0)=\Phi_{ext,m}(\cdot\,,0)$,}
\\
\\
\displaystyle
\mbox{\boldmath $v$}_{int,m}, & \text{if $\Phi(\,\cdot\,,0)=\Phi_{int,m}(\cdot\,,0)$}.
\end{array}
\right.
$$
In what follows we write
$$\displaystyle
\tilde{\tau}=\sqrt{\tau}\cdot\sqrt{\frac{\rho}{\mu}}.
$$
We have the following expression.
\proclaim{\noindent Proposition 3.1.}
(i)  Let $\vert x-z\vert>\eta$.  We have
$$\begin{array}{ll}
\displaystyle
\mbox{\boldmath $v$}_{ext,m}(x)
&
\displaystyle
=\frac{\eta^{2(1+m)}}{\mu}
\cdot\frac{e^{-\tilde{\tau}\,\vert x-z\vert}}{\,\vert x-z\vert}\,
a_m(\tilde{\tau})
\,\mbox{\boldmath $a$},
\end{array}
$$
where
$$\displaystyle
a_m(\tilde{\tau})
=\frac{1}{\tilde{\tau}}\int_0^1s(1-s^2)^m\sinh(\eta\tilde{\tau}\,s)\,ds.
$$

(ii)  Let $\vert x-z\vert<R_1$.  We have
$$\displaystyle
\mbox{\boldmath $v$}_{int,m}(x)
=\frac{2}{\mu}
\cdot\frac{\sinh\,\tilde{\tau}\,\vert x-z\vert}{\vert x-z\vert}
\,b_m(\tilde{\tau})\,\mbox{\boldmath $a$},
$$
where
$$\displaystyle
b_m(\tilde{\tau})
=\frac{1}{\tilde{\tau}}\,\int_{R_1}^{R_2}\,s(R_2^2-s^2)^m(s^2-R_1^2)^m\,e^{-s\tilde{\tau}}\,ds.
$$

\endproclaim

The proof is given in the next section.

Concerning with the asymptotic behaviour of the coefficients $a_m(\tilde{\tau})$
and $b_m(\tilde{\tau})$, by Theorem 7.1 on p.81 in \cite{Ol}, we have,
as $\tilde{\tau}\rightarrow\infty$
$$
\displaystyle
a_m(\tilde{\tau})
\sim
\eta\,2^{m-1}m!\frac{e^{\tilde{\tau}\,\eta}}{(\tilde{\tau}\eta)^{m+2}}
\tag {3.8}
$$
and
$$
\\
\displaystyle
b_m(\tilde{\tau})
\sim
2^m\,m!R_1^{m+1}\,(R_2^2-R_1^2)^m\frac{e^{-R_1\tilde{\tau}}}{\tilde{\tau}^{m+2}}.
\tag{3.9}
$$

Now we are ready to prove the key lemma stated below.

\proclaim{\noindent Lemma 3.2.}
There exist positive numbers $\tau_0$ and $\kappa$ such that, for all $\tau\ge\tau_0$
$$\displaystyle
\tau^{\kappa}e^{2\sqrt{\tau}\sqrt{\rho/\mu}\,\text{dist}\,(D,K)}\,J(\tau)\ge C.
\tag {3.10}
$$

\endproclaim
{\it\noindent Proof.}
It follows from (3.7) that
$$\displaystyle
J(\tau)\ge\frac{1}{2}\tau\rho\Vert \mbox{\boldmath $w$}_{00}\Vert_{L^2(D)}^2
+O(\tau^{-1}e^{-2\tau T}).
\tag {3.11}
$$
Thus it suffices to give a lower estimate for $\Vert\mbox{\boldmath $w$}_{00}\Vert_{L^2(D)}$.

First consider the case when $K=\overline B_{\eta}$.
It follows from (3.4) and (i) in Proposition 3.1 that $\mbox{\boldmath $w$}_{00}$ takes the form
$$\begin{array}{ll}
\displaystyle
\mbox{\boldmath $w$}_{00}(x)=
\frac{\eta^{2(1+m)}}
{\mu}
\cdot a_m(\tilde{\tau})
\nabla
\left(
\frac{e^{-\tilde{\tau}\,\vert x-z\vert}}{\vert x-z\vert}
\right)
\times\mbox{\boldmath $a$}, & x\in D.
\end{array}
\tag {3.12}
$$
Here we have
$$\displaystyle
\nabla
\left(
\frac{e^{-\tilde{\tau}\,\vert x-z\vert}}{\vert x-z\vert}
\right)
=-e^{-\tilde{\tau}\,\vert x-z\vert}
\left(\frac{1}{\vert x-z\vert^2}+\frac{1}{\vert x-z\vert}\right)\frac{x-z}{\vert x-z\vert}.
$$
Thus one gets
$$\displaystyle
\int_D
\left\vert\nabla
\left(
\frac{e^{-\tilde{\tau}\,\vert x-z\vert}}{\vert x-z\vert}
\right)
\times\mbox{\boldmath $a$}
\right\vert^2
\,dx
\ge C^2
\int_De^{-2\tilde{\tau}\,\vert x-z\vert}
\left\vert\frac{x-z}{\vert x-z\vert}\times\mbox{\boldmath $a$}\right\vert^2\,dx,
$$
where
$$\displaystyle
C=\frac{1}{\sup_{x\in D}\,\vert x-z\vert^2}+\frac{1}{\sup_{x\in D}\,\vert x-z\vert}.
$$
By Lemma A.3 in \cite{Ithermo} there exist positive constants $C_1$, $\tau_0$ and $\kappa'$ such that
$$\displaystyle
\int_D
e^{-2\tilde{\tau}\,\vert x-z\vert}
\left\vert\frac{x-z}{\vert x-z\vert}\times\mbox{\boldmath $a$}\right\vert^2\,dx
\ge C_1e^{-2\eta\tilde{\tau}}\cdot\tau^{-\kappa'}e^{-2\tilde{\tau}\,\text{dist}\,(D,K)}
$$
for all $\tau\ge\tau_0$.
Thus from (3.12) one gets
$$\displaystyle
\Vert\mbox{\boldmath $w$}_{00}\Vert_{L^2(D)}^2\ge
\left\{\frac{\eta^{2(1+m)}}
{\mu}
\cdot a_m(\tilde{\tau})\right\}^2
C^2C_1
\,e^{-2\eta\tilde{\tau}}\cdot\tau^{-\kappa'}e^{-2\tilde{\tau}\,\text{dist}\,(D,K)}.
$$
Now from this together with (3.8) and (3.11) we see that (3.10) is valid by choosing $\kappa=\frac{m+2}{2}+\kappa'-1$.

Next consider the case when $K=\overline{B_{R_2}}\setminus B_{R_1}$.
From (3.4) and (ii) in Proposition 3.1
we have 
$$\begin{array}{ll}
\displaystyle
\mbox{\boldmath $w$}_{00}(x)
=\frac{2}{\mu}\cdot b_m(\tilde{\tau})
\nabla
\left(
\frac{\sinh\,\tilde{\tau}\,\vert x-z\vert}{\vert x-z\vert}
\right)
\times
\mbox{\boldmath $a$}, & x\in D.
\end{array}
\tag {3.13}
$$
By the proof of Lemma 4.3 in \cite{IE06}
we know that there exist positive constants $C_2$, $\tau_1$ and $\kappa''$ such that
$$\displaystyle
\int_D
\left\vert
\nabla
\left(
\frac{\sinh\,\tilde{\tau}\,\vert x-z\vert}{\vert x-z\vert}
\right)
\times
\mbox{\boldmath $a$}
\right\vert^2\,dx
\ge C_2
\tau^{-\kappa''}\,e^{2\tilde{\tau}\,R_D(z)}
$$
for all $\tau\ge\tau_1$.
This together with (3.13) yields
$$\displaystyle
\Vert\mbox{\boldmath $w$}_{00}\Vert_{L^2(D)}^2
\ge
\left\{\frac{2}{\mu}\cdot b_m(\tilde{\tau})\right\}^2
C_2\tau^{-\kappa''}\,e^{2\tilde{\tau}\,R_D(z)}.
$$
Now from this together with (3.9) and (3.11) we see that (3.10) is valid by choosing 
$\kappa=\frac{m+2}{2}+\kappa''-1$.  Note that we have made use of the relationship
$\text{dist}\,(D,K)=R_1-R_D(z)$.

\noindent
$\Box$

$\quad$

{\bf\noindent Remark 3.1.}
In contrast to the wave equation case, for example, (ii) of Lemma 2.4 in \cite{Iwave}, there is no restriction on the lower bound for $T$ in Lemma 3.2.
See $J(\tau)$ on (2.6), which comes from $\mbox{\boldmath $w$}_0$.  The $\mbox{\boldmath $w$}_0$ comes
from $\mbox{\boldmath $v$}$ which solves (1.3) and is generated by the initial data supported on $K$.  
The $K$ is remote from $D$.
So Lemma 3.2 implicitly indicates that the signal governed by the Stokes system propagates with infinity speed.

\subsection{Finishing the proof of Theorem 1.1}

First we prepare three estimates.
Letting $U=\Omega$ in Lemma 3.1, we have
$$\displaystyle
\Vert\mbox{\boldmath $w$}_0\Vert_{L^2(\Omega)}=O(\sqrt{\tau}\,e^{-\sqrt{\tau}\,\sqrt{\rho/\mu}\,\text{dist}\,(\Omega,K)}
+\tau^{-1}e^{-\tau T}).
\tag {3.14}
$$
Besides a combination of (3.1) in the case $U=D$ and (2.6) yields
$$\displaystyle
J(\tau)=O(\tau^2 e^{-2\sqrt{\tau}\sqrt{\rho/\mu}\,\text{dist}\,(D,K)}+\tau^{-1}e^{-2\tau T}).
\tag {3.15}
$$
From this together with (2.13) we obtain
$$\displaystyle
E(\tau)=O(\tau^3 e^{-2\sqrt{\tau}\sqrt{\rho/\mu}\,\text{dist}\,(D,K)}+e^{-2\tau T}).
\tag {3.16}
$$

Applying (3.14) and (3.15) to the right-hand side on (2.14), we obtain
$$\displaystyle
e^{-\tau T}{\cal R}(\tau)=O(e^{-2^{-1}\tau T}).
$$
Note that we just used (3.14) in which $\text{dist}\,(\Omega,K)$ is replaced with $0$.
Then (2.5) gives
$$\displaystyle
I_T(\tau;\mbox{\boldmath $v$}, 0)
=J(\tau)+E(\tau)+O(e^{-2^{-1}\tau T}).
\tag {3.17}
$$

Since $E(\tau)\ge 0$, it follows from (3.10) that there exists a positive 
number $\tau_0$ such that, for all $\tau\ge\tau_0$
$$\displaystyle
\tau^{\kappa}e^{2\sqrt{\tau}\,\sqrt{\rho/\mu}\,\text{dist}\,(D,K)}\,I_T(\tau;\mbox{\boldmath $v$}, 0)
\ge C/2.
\tag {3.18}
$$

Applying (3.15) and (3.16) to the right-hand side on (3.17), we obtain
$$\displaystyle
e^{2\sqrt{\tau}\,\sqrt{\rho/\mu}\,\text{dist}\,(D,K)}\,I_T(\tau;\mbox{\boldmath $v$}, 0)
=O(\tau^3).
\tag {3.19}
$$
Then, from (3.18) and (3.19) we immediately obtain (1.8) and (1.9).

This completes the proof of Theorem 1.1.

\section{Proof of Proposition 3.1}

The computation presented here is a combination of the Parseval identity and the residue calculus as done in the proof
of Proposition 3.1 in \cite{IHeatWave}.

Write
$$\begin{array}{ll}
\displaystyle
\int_{\Bbb R^3}\,e^{-ix\cdot\xi}\,\Phi_{ext,m}(x)\,dx
&
\displaystyle
=\int_{B_{\eta}}e^{-ix\cdot\xi}\,(\eta^2-\vert x-z\vert^2)^m\,dx
\\
\\
\displaystyle
&
\displaystyle
=e^{-iz\cdot\xi}\,
\int_0^{\eta}r^2(\eta^2-r^2)^m\,dr
\int_{S^2}e^{-ir\omega\cdot\xi}\,d\omega.
\end{array}
$$
Since we have
$$\displaystyle
\int_{S^2}\,e^{-ir\omega\cdot\xi}\,d\omega=4\pi\cdot\frac{\sin r\vert\xi\vert}{r\vert\xi\vert},
$$
one gets
$$\begin{array}{ll}
\displaystyle
\hat{\Phi_{ext,m}}(\xi)
&
\displaystyle
=4\pi\,e^{-iz\cdot\xi}
\,\int_0^{\eta}r^2(\eta^2-r^2)^m\frac{\sin r\vert\xi\vert}{r\vert\xi\vert}\,dr
\\
\\
\displaystyle
&
\displaystyle
=4\pi\eta^{2(1+m)}\,e^{-iz\cdot\xi}\vert\xi\vert^{-1}
A(\vert\xi\vert),
\end{array}
\tag {4.1}
$$
where
$$\displaystyle
A_m(\zeta)=\int_0^1s(1-s^2)^m\sin\,(\eta\,\zeta\,s)\,ds.
$$
Similary we have
$$\begin{array}{ll}
\displaystyle
\hat{\Phi_{int,m}}(\xi)
&
\displaystyle
=4\pi\,e^{-iz\cdot\xi}\,\vert\xi\vert^{-1}
B_m(\vert\xi\vert),
\end{array}
\tag {4.2}
$$
where
$$\displaystyle
B_m(\zeta)
=
\int_{R_1}^{R_2}s(R_2^2-s^2)^m(s^2-R_1^2)^m\,
\sin\,\zeta s\,ds.
$$

Using the Parseval identity, (3.3) and (4.1) we obtain
$$\begin{array}{ll}
\displaystyle
\mbox{\boldmath $v$}_{ext,m}(x)
&
\displaystyle
=\frac{1}{\mu(2\pi)^3}
\int_{\Bbb R^3}\,\frac{e^{ix\cdot\xi}}{\vert\xi\vert^2+\tilde{\tau}^2}
\cdot
4\pi\eta^{2(1+m)}\,e^{-iz\cdot\xi}\vert\xi\vert^{-1}
A_m(\vert\xi\vert)
\,d\xi\,\mbox{\boldmath $a$}
\\
\\
\displaystyle
&
\displaystyle
=\frac{2\eta^{2(1+m)}}{\mu(2\pi)^2}
\int_{\Bbb R^3}\,\frac{e^{i\,\vert x-z\vert\,\xi_3}}{\vert\xi\vert(\vert\xi\vert^2+\tilde{\tau}^2)}\,
A_m(\vert\xi\vert)
\,d\xi\,\mbox{\boldmath $a$}
\\
\\
\displaystyle
&
\displaystyle
=\frac{\eta^{2(1+m)}}{\mu\pi}
\int_0^{\infty}\,\frac{r}{r^2+\tilde{\tau}^2}\,
A_m(r)\,
dr
\,\int_0^{\pi}\sin\phi\,
e^{i\,\vert x-z\vert\,r\cos\phi}
\,d\phi
\,
\mbox{\boldmath $a$}
\\
\\
\displaystyle
&
\displaystyle
=\frac{\eta^{2(1+m)}}{\mu\pi}
\cdot\frac{1}{i\vert x-z\vert}
\int_0^{\infty}\,\frac{1}{r^2+\tilde{\tau}^2}\,
A_m(r)\,
dr
\,(e^{i\vert x-z\vert\,r}-e^{-i\vert x-z\vert\,r})
\,
\mbox{\boldmath $a$}
\\
\\
\displaystyle
&
\displaystyle
=\frac{2\eta^{2(1+m)}}{\mu\pi}
\cdot\frac{1}{\vert x-z\vert}
\int_0^{\infty}\,\frac{1}{r^2+\tilde{\tau}^2}\,
A_m(r)\,
\sin\,\vert x-z\vert\,r
\,dr
\,
\mbox{\boldmath $a$}
\\
\\
\displaystyle
&
\displaystyle
=\frac{\eta^{2(1+m)}}{\mu\pi}
\cdot\frac{1}{\vert x-z\vert}
\int_{-\infty}^{\infty}
\,\frac{1}{r^2+\tilde{\tau}^2}\,
A_m(r)\,
\sin\,\vert x-z\vert\,r
\,dr
\,
\mbox{\boldmath $a$}.
\end{array}
\tag {4.3}
$$
Note that, at the final step we made use of the evenness of the integrand with respect to $r$.

Let $\zeta=Re^{i\theta}$, $0\le\theta\le\pi$.
We have
$$\displaystyle
\left\vert\,A_m(\zeta)\right\vert
\le Ce^{\eta\,R\sin\theta}
$$
and
$$\displaystyle
\vert\,e^{i\vert x-z\vert\,\zeta}\vert=e^{-\vert x-z\vert R\sin\theta}.
$$
A standard residue calculus yields:  if $\vert x-z\vert>\eta$, then
we have
$$\begin{array}{l}
\displaystyle
\,\,\,\,\,\,
\int_{-\infty}^{\infty}
\,\frac{1}{r^2+\tilde{\tau}^2}\,
A_m(r)\,
e^{i\vert x-z\vert r}
\,dr
\\
\\
\displaystyle
=2\pi i
\,\frac{1}{2i\tilde{\tau}}\,
A_m(i\tilde{\tau})\,
e^{i\vert x-z\vert i\tilde{\tau}}
\\
\\
\displaystyle
=
\,\frac{\pi}{\tilde{\tau}}\,
\int_0^1s(1-s^2)^m\sin\,\eta\,i\tilde{\tau}\,s\,ds\,
e^{-\vert x-z\vert\tilde{\tau}}
\\
\\
\displaystyle
=i\,\frac{\pi}{\tilde{\tau}}\,
\int_0^1s(1-s^2)^m\sinh(\eta\tilde{\tau}\,s)\,ds\,e^{-\vert x-z\vert\tilde{\tau}}
\end{array}
$$
and thus
$$
\displaystyle
\int_{-\infty}^{\infty}
\,\frac{1}{r^2+\tilde{\tau}^2}\,
A_m(r)
\,
\sin\,\vert x-z\vert\,r
\,dr
=\pi a_m(\tilde{\tau})\,\,e^{-\vert x-z\vert\tilde{\tau}}.
$$
Now from (4.3) we obtain (i).

For (ii), using the Parseval identity, (3.3) and (4.2) we obtain
$$\begin{array}{ll}
\displaystyle
\mbox{\boldmath $v$}_{int,m}(x)
&
\displaystyle
=\frac{1}{\mu(2\pi)^3}
\int_{\Bbb R^3}\,\frac{e^{ix\cdot\xi}}{\vert\xi\vert^2+\tilde{\tau}^2}
\cdot
4\pi\,\vert\xi\vert^{-1}\,e^{-iz\cdot\xi}\,
B_m(\vert\xi\vert)
d\xi\,
\mbox{\boldmath $a$}
\\
\\
\displaystyle
&
\displaystyle
=\frac{4\pi}{\mu(2\pi)^3}
\int_{\Bbb R^3}\,\frac{e^{i\vert x-z\vert\,\xi_3}}{\vert\xi\vert(\vert\xi\vert^2+\tilde{\tau}^2)}\,
B_m(\vert\xi\vert)
d\xi
\,\mbox{\boldmath $a$}
\\
\\
\displaystyle
&
\displaystyle
=\frac{2}{\mu\pi}
\int_0^{\infty}\frac{r}{r^2+\tilde{\tau}^2}\cdot B_m(r)\,dr
\int_0^{\pi}\,
e^{i\vert x-z\vert\,r\cos\phi}\,\sin\phi\,d\phi
\,\mbox{\boldmath $a$}
\\
\\
\displaystyle
&
\displaystyle
=-i\frac{2}{\mu\pi}\frac{1}{\vert x-z\vert}
\int_0^{\infty}\frac{1}{r^2+\tau^2}\cdot B_m(r)
(e^{i\vert x-z\vert\,r}-e^{-i\vert x-z\vert\,r})\,dr
\,
\mbox{\boldmath $a$}
\\
\\
\displaystyle
&
\displaystyle
=\frac{4}{\mu\pi}\frac{1}{\vert x-z\vert}
\int_0^{\infty}\frac{1}{r^2+\tau^2}\cdot B_m(r)
\sin\,\vert x-z\vert r\,dr
\,
\mbox{\boldmath $a$}
\\
\\
\displaystyle
&
\displaystyle
=\frac{2}{\mu\pi}\frac{1}{\vert x-z\vert}
\int_{-\infty}^{\infty}\frac{1}{r^2+\tau^2}\cdot B_m(r)
\sin\,\vert x-z\vert r\,dr
\,
\mbox{\boldmath $a$}.
\end{array}
\tag {4.4}
$$
Here one can write
$$\displaystyle
\int_{-\infty}^{\infty}\frac{1}{r^2+\tau^2}\cdot B_m(r)
\sin\,\vert x-z\vert r\,dr
=
\text{Im}\,
\int_{-\infty}^{\infty}\frac{C_m(r)}{r^2+\tau^2}
\sin\,\vert x-z\vert r\,dr,
$$
where
$$\begin{array}{ll}
\displaystyle
C_m(\zeta)
=\int_{R_1}^{R_2}s(R_2^2-s^2)^m(s^2-R_1^2)^m\,
e^{is\zeta}\,ds,
& m=0,1,\cdots.
\end{array}
$$

Let $\zeta=Re^{i\theta}$, $0\le\theta\le\pi$.
We have
$$\displaystyle
\vert C_m(z)\vert\le C(R_1,R_2)e^{-RR_1\sin\theta}
$$
and
$$\displaystyle
\vert\sin\vert x-z\vert\zeta\vert\le e^{R\vert x-z\vert\sin\theta}.
$$
Using the residue theorem, we obtain: if $\vert x-z\vert<R_1$, then
$$\displaystyle
\int_{-\infty}^{\infty}
\frac{C_m(r)}{r^2+\tilde{\tau}^2}
\sin\,\vert x-z\vert r\,dr
=\pi\tilde{\tau}^{-1}C_m(i\tilde{\tau})\sin\,(\vert x-z\vert\,i\tilde{\tau}).
$$
Note
$$\displaystyle
C_m(i\tilde{\tau})=\int_{R_1}^{R_2}s(R_2^2-s^2)^m(s^2-R_1^2)^m\,e^{-s\tilde{\tau}}\,ds
$$
and
$$\displaystyle
\sin\,(\vert x-z\vert\,i\tilde{\tau})
=\frac{e^{-\tilde{\tau}\,\vert x-z\vert}-e^{\tilde{\tau}\vert x-z\vert}}{2i}
=i\,\frac{1}{2}(e^{\tilde{\tau}\vert x-z\vert}-e^{-\tilde{\tau}\vert x-z\vert})
=i\sinh\,\tilde{\tau}\vert x-z\vert.
$$
Thus one gets
$$\displaystyle
\text{Im}\,
\int_{-\infty}^{\infty}
\frac{C_m(r)}{r^2+\tilde{\tau}^2}
\sin\,\vert x-z\vert r\,dr
=\pi\,b_m(\tilde{\tau})
\,\sinh\,\tilde{\tau}\vert x-z\vert.
$$
Now from (4.4) we obtain (ii).

This completes the proof of Proposition 3.1.

$\quad$

\centerline{{\bf Acknowledgments}}

The author was partially supported by Grant-in-Aid for
Scientific Research (C)(No. 17K05331) and (B)(No. 18H01126) 
of Japan  Society for
the Promotion of Science.

$\quad$

\end{document}